\documentclass[a4paper,12pt]{article}

\usepackage[latin1]{inputenc}
\usepackage{abstract} 
\usepackage[arrow, matrix, curve]{xy}
\usepackage{amsmath,amssymb,amsfonts,amsthm,amsopn}
\usepackage{nicefrac}
\usepackage{fancyhdr}

\title{A Dynamical Bogomolov Property}
\author{Lukas Pottmeyer}
\date{\today}

\DeclareMathOperator{\val}{val}

\DeclareMathOperator{\supp}{supp}

\DeclareMathOperator{\Aut}{Aut}

\DeclareMathOperator{\PrePer}{PrePer}

\newcommand{\Monat}{%
\ifcase\month
 Monat 0 \or Januar \or Februar \or März  \or April \or Mai \or Juni \or Juli%
 \or August \or September \or Oktober \or November \or Dezember
\fi}

\newtheorem{Theorem}{Theorem}[section]

\newtheorem{Proposition}[Theorem]{Proposition}
\newtheorem{Korollar}[Theorem]{Corollary}

\newtheorem{Definition}[Theorem]{Definition}
\newtheorem{Bemerkung}[Theorem]{Remark}

\makeindex

\begin{document}

\maketitle

\begin{abstract}
A field $F$ is said to have the Bogomolov Property related to a height function $h$ if $h(\alpha)$ is either $0$ or bounded from below by a positive constant for all $\alpha \in F$. In this paper we prove that the maximal algebraic extension of a number field $K$, which is unramified at a place $v\mid p$, has the Bogomolov Property related to all canonical heights coming from a Lattès map related to a Tate elliptic curve. To prove this algebraical statement we use analytic methods on the related Berkovich spaces.
\end{abstract}

In the whole paper $h$ is the standard logarithmic height, $p$ a prime number, $K$ a number field and $K^{nr,v}$ the maximal algebraic extension of $K$, whicht is unramified at the place $v\mid p$. 

\begin{Definition}
Let $\varphi: \mathbb{P}^{n} \rightarrow \mathbb{P}^{n}$ be a morphism of degree $d\geq2$. The canonical height associated to $\varphi$ is the unique height function $$\widehat{h}_{\varphi}:\mathbb{P}^{n} \rightarrow \mathbb{R}$$
with the properties

$$\begin{array}{lcr} 
\widehat{h}_{\varphi}(P)=h(P) + O(1) & \text{ and } & \widehat{h}_{\varphi}(\varphi(P))=d\widehat{h}_{\varphi}(P).
\end{array}$$

\end{Definition}

See [Si07], 3.20, for a proof of the existence and uniquness of this height. In [Si07], 3.22, it is shown that
$$\widehat{h}_{\varphi}(P)=0 \Leftrightarrow P \text{ is a preperiodic point of } \varphi .$$
We will write $\PrePer(\varphi)$ for the set of all preperiodic points of $\varphi$.

\begin{Definition}
Let $E$ be an elliptic curve over a field $F$ of characteristic $0$, $\Gamma$ a non trivial subgroup of $\Aut(E)$, $\varphi$ an endomorphism of $E$ and $\pi: E \rightarrow \nicefrac{E}{\Gamma}\cong\mathbb{P}_{F}^{1}$ a finite covering. A morphism $f:\mathbb{P}_{F}^{1} \rightarrow \mathbb{P}_{F}^{1}$ is called Lattès map related to $E$, if we have a commutative diagram
$$\begin{xy}
  \xymatrix{
      E \ar[r]^\varphi \ar[d]_\pi    &   E \ar[d]^\pi  \\
      \nicefrac{E}{\Gamma} \ar[r]^\varphi \ar[d]_T &  \nicefrac{E}{\Gamma} \ar[d]^T \\
      \mathbb{P}_{F}^{1} \ar[r]^f             &   \mathbb{P}_{F}^{1}   
  }
\end{xy}$$
\end{Definition}
Here $T$ is an isomorphism of $\nicefrac{E}{\Gamma}$ and $\mathbb{P}_{F}^{1}$.

\begin{Bemerkung} \rm
This definition is independent of the choice of the isomprphism $T$. A change of the isomorphism is equivalent to a change of the coordinate on $\mathbb{P}_{F}^{1}$. So it would just change the representation of the map $f$ but not the map itself. This allows us to consider just the reduced Lattès diagram
$$\begin{xy}
  \xymatrix{
      E \ar[r]^\varphi \ar[d]_\pi    &   E \ar[d]^\pi  \\
      \mathbb{P}_{F}^{1} \ar[r]^f             &   \mathbb{P}_{F}^{1}   
  }
\end{xy}$$
with a finite covering $\pi:E \rightarrow \mathbb{P}_{F}^{1}$.

\end{Bemerkung}

If $f$ is a Lattès map for a Tate elliptic curve $E$, then we will prove that $K^{nr,v}$ has the Bogomolov Property related to $\widehat{h}_f$. In a slightly different context the Bogomolov Property for $K^{nr,v}$ is already known. We have the following result.

\begin{Theorem}[Gubler]
Let $A$ be an abelian variety over $K$, which is totally degenerate at $v$. Further let $L$ be an ample even line bundle and $K'$ be a finite extension of $K^{nr,v}$. Then there is a $\varepsilon >0$, such that $\widehat{h}_{L}(P) \geq \varepsilon$ for all non-torsion points $P \in A(K')$.
\end{Theorem}
\textbf{Proof:} See [Gu07], Corollary 6.7. $\hfill \square$

\vspace{0.4cm}

In this theorem $\widehat{h}_L $ stands for the canonical height, also called Néron-Tate height, associated to the line bundle $L$. For details on this height we refer to [BG], Chapter 9.2, and [La], Chapter 5.

One might expect the Bogomolov Property of $K^{nr,v}$ related to $\widehat{h}_f$, for a Lattès map $f$ associated to a Tate elliptic curve, to be a direct consequence of theorem 1.4.

Let $A$ be a Tate elliptic curve $y^2 +xy = x^3 + a_4 x +a_6$ over $K$, which is totally degenerate at $v$. Notice that the Tate elliptic curves are precisely those with total degeneration at $v$. Let $\pi$ be the projection on the $x$-coordinate, $\varphi = [m]$, $m \in \mathbb{Z}$, be the multiplication by $m$ map and $f$ the corresponding Lattès map. (Later we will show that this is the general case.) The multiplication by $[-1]$ on $A$ is given by $[-1](x,y)=(x,-x-y)$ (for example [BG], Proposition 8.3.8), so we can deduce that the line bundle $L:=\pi^* \mathcal{O}(1)$ is even. By 1.1, 1.2 and the universal property of $\widehat{h}_L$ we get the equation 
\begin{eqnarray} \widehat{h}_f \circ \pi = \widehat{h}_L . \label{height} \end{eqnarray}
Now we see that the Bogomolov Property of $K^{nr,v}$ related to $\widehat{h}_f$ would be a direct consequence of Theorem 1.4 if $[K^{nr,v}(\pi^{-1}(\mathbb{P}_{K}^{1}(K^{nr,v}))):K^{nr,v}]$ is finite.
But in general this is not the case. A counterexample for $K=\mathbb{Q}$ is the Tate elliptic curve $y^2 +xy = x^3 +6$ when we choose $p=3$. To see this one needs a lot of algebraic number theory, but the proof will not be shown here.
After this short justification we will prove our main theorem.

\begin{Theorem}
Let $f$ be a Lattès map related to a Tate elliptic curve over $K^{nr,v}$ and let $\widehat{h}_f$ be the canonical height associated to $f$. Then $K^{nr,v}$ has the Bogomolov Property related to $\widehat{h}_f$. In other words:
\begin{center}
There is a $C>0$, such that for all $P \in K^{nr,v} \setminus \PrePer(f)$ we have $\widehat{h}_f (P)>C$.
\end{center}
Moreover there are just finitely many points $a$ in $K^{nr,v}$ with $\widehat{h}_f (a)=0$.
\end{Theorem}

Let $w \mid v$  be an absulute value on $K^{nr,v}$. By $\mathbb{K}$ we denote the completion of the algebraic closure of the completion of $K^{nr,v}$ by $w$. We denote the unique extension of $w$ to $\mathbb{K}$ also by $w$. Notice that $\mathbb{K}$ is algebraically closed ([BGR], Proposition 3.4.3). Let $E : y^2 + xy = x^3 + a_4 x + a_6$ be a Tate elliptic curve over $\mathbb{K}$ and $\pi: E \rightarrow \mathbb{P}_{\mathbb{K}}^{1}$ be a finite covering. Further we take a Lattès map $f$ with the commutative diagram

$$\begin{xy}
  \xymatrix{
      E \ar[r]^\varphi \ar[d]_\pi    &   E \ar[d]^\pi  \\
      \mathbb{P}^{1}_{\mathbb{K}} \ar[r]^f             &   \mathbb{P}^{1}_{\mathbb{K}}   
  }
\end{xy}$$

where $\varphi$ is the multiplication by $m\in\mathbb{Z}$ with $\vert m \vert \geq 2$. The GAGA-functor on Berkovich spaces (see [Ber], Chapter 3.4) leads us to the commutative diagram

$$\begin{xy}
  \xymatrix{
      E^{an} \ar[r]^{\varphi^{an}} \ar[d]_{\pi^{an}}    &   E^{an} \ar[d]^{\pi^{an}}  \\
      (\mathbb{P}^{1}_{\mathbb{K}})^{an} \ar[r]^{f^{an}}             &   (\mathbb{P}^{1}_{\mathbb{K}})^{an}   
  }
\end{xy}$$

By $E^{an}$ and $(\mathbb{P}_{\mathbb{K}}^{1})^{an}$ we denote the Berkovich spaces related to $E$, respectevely $\mathbb{P}_{\mathbb{K}}^{1}$, and by $\varphi^{an}$, $\pi^{an}$ and $f^{an}$ we denote the analytification (in the sence of Berkovich) of the respective map.

We will work with the valuation function on $(\mathbb{P}_{\mathbb{K}}^{1})^{an}$

$$ \begin{array}{ccc}
\val : (\mathbb{P}^{1}_{\mathbb{K}})^{an} \rightarrow \mathbb{R}\cup\{\pm \infty\} & ; & y \mapsto  -\log \vert X \vert_{y},
   \end{array}$$

where $X$ is the variable in the ring of polynomials $\mathbb{K}[X]$. As an analytic group $E$ is isomorphic to $\nicefrac{\mathbb{K}^*}{q^{\mathbb{Z}}}$ for a $q \in \mathbb{K}$ with $\vert q \vert_w < 1$. So there is a canonical valuation function $\overline{\val}$ on $E^{an}=\nicefrac{(\mathbb{G}^{1}_{m})^{an}}{q^\mathbb{Z}}$

$$ \begin{array}{ccc}
\overline{\val} : \nicefrac{(\mathbb{G}^{1}_{m})^{an}}{q^\mathbb{Z}} \rightarrow \nicefrac{\mathbb{R}}{w(q) \mathbb{Z}} & ; & \overline{y} \mapsto  \overline{-\log \vert X \vert_{y}}.
   \end{array}$$
   
Obviously we have $\overline{\val} (E^{an})=\nicefrac{\mathbb{R}}{w(q)\mathbb{Z}}$.

\begin{Proposition}
$E$ has no complex multiplication and the map $\pi: E \rightarrow \mathbb{P}_{\mathbb{K}}^{1}$ is, after a suitable coordinate transformation, explicitly given by $(x,y) \mapsto x$.
\end{Proposition}

\textbf{Proof:} Let $j(E)$ be the $j$-invariant of $E$. Then we have the equation
$$\vert j(E)\vert_w = \vert q \vert_{w}^{-1} > 1.$$
Since $w$ is an extension of a $p$-adic absolute value, $j(E)$ is no algebraic integer. That shows the first statement (see [Si99], Theorem II.6.1). Further we know $\Aut(E)=\{id,[-1]\}$ (see [Si07], Proposition 6.26). So there is just one quotient curve of $E$ which is different from $E$ itself, namely $\nicefrac{E}{\Aut(E)}\cong\mathbb{P}_{\mathbb{K}}^{1}$. The function field of $\nicefrac{E}{\Aut(E)}$ is as the fixed field $\mathbb{K}(x,y)^{\Aut(E)}$ given by $\mathbb{K}(x,y^2 + xy)=\mathbb{K}(x)$. Therefore the projection $\pi:E\rightarrow\nicefrac{E}{\Aut(E)}$ is given by $(x,y)\mapsto x$. \hfill $\square$

\vspace{0.4cm}

Again we take the even ample line bundle $L:=\pi^* \mathcal{O}(1)$. We have $\deg(f)=\deg(\varphi)=m^2$ (see [Si07], Theorem 6.51). This leads us to an isomorphism $\Phi: \mathcal{O}(1)^{m^2} \rightarrow f^* \mathcal{O}(1)$. $L$ is even and so the theorem of the cube tells us $\varphi^* L \cong L^{m^2}$. We choose the isomorphism $\Psi:=\pi^* \Phi$. Notice that, by the commutativity of the Lattès diagram, we have $\varphi^* L=\varphi^* \pi^* \mathcal{O}(1)=\pi^* f^* \mathcal{O}(1)$. There are unique metrics $\|.\|_f$ and $\|.\|_\varphi$ on $\mathcal{O}(1)$ respectively $L$ (to be more formal: on the analytifications of these line bundles) with the properties
$$(\Phi^{an})^* (f^{an})^*\|.\|_{f}=\|.\|^{m²}_{f} \text{ and } (\Psi^{an})^* (\varphi^{an})^*\|.\|_{\varphi}=\|.\|_{\varphi}^{m²}.$$
As usual we denote by $\Psi^{an}$ and $\Phi^{an}$ the analytifications (in the sence of Berkovich) of $\Psi$ and $\Phi$. For the existence of these metrics and further information we refer to [Zh], [Gu98] (especially Theorem 7.12) and [Gu10]. Just using the definitions of the different maps we get
$$(\Psi^{an})^{*}(\varphi^{an})^* (\pi^{an})^* \|.\|_{f} =((\pi^{an})^* \|.\|_{f})^{m²} .$$
This implies the equation \begin{eqnarray}
 (\pi^{an})^* \|.\|_{f} = \| . \|_{\varphi} . \label{metrics} 
\end{eqnarray}

Now we have a look at the canonical measures, also called Chambert-Loir measures, of the arithmetical dynamical systems $(E,\varphi,L)$ and $(\mathbb{P}_{\mathbb{K}}^{1},f,\mathcal{O}(1))$. We denote these measures by $\mu_\varphi = c_{1}(L,\| . \|_{\varphi})$ and $\mu_f = c_{1}(\mathcal{O}(1),\| . \|_{f})$. For details, we refer to [Ch] and [Gu10]. Using the projection formula (for example [Gu07a], Corollary 3.9 b)) and \eqref{metrics} we deduce
$$(\pi^{an})_* \mu_{\varphi} = \deg(\pi) \mu_f .$$
This leads us to $ \supp((\pi^{an})_* \mu_{\varphi})=\supp(\mu_f )$. As $\mu_\varphi$ is a positive measure we get
\begin{eqnarray}
 \pi^{an}(\supp(\mu_\varphi))=\supp((\pi^{an})_* \mu_{\varphi})=\supp(\mu_f ) .\label{support}
\end{eqnarray}

\begin{Bemerkung} \rm
Every disk $(a,r)$ around $a \in \mathbb{K}$ with radius $r\in \mathbb{R}^+$ gives us a multiplicative seminorm on the Tate-algebra $\mathbb{K}\{X,qX^{-1}\}$, and hence a point of $E^{an}$. Explicitly $\vert.\vert_{(a,r)}$ is given by
$$\vert \sum_{n\in \mathbb{Z}} a_n X^n \vert_{(a,r)} = \vert \sum_{n\in \mathbb{Z}} b_n (X-a)^n \vert_{(a,r)} = \max_{n\in \mathbb{Z}} \vert b_n \vert_w r^n .$$
\end{Bemerkung}

The subdomain of  $E^{an}$ consisting of all points $(0,r)$, with $\vert q \vert_w < r < 1$ is called the skeleton of $E$. We denote the skeleton of $E$ by $S(E)$. It is easy to see, that $\overline{\val}$ maps $S(E)$ homeomorphic onto $\nicefrac{\mathbb{R}}{w(q)\mathbb{Z}}$. For the general theory of skeletons we refer to [Ber], Chapter 6.5, and for more information on our special case we refer to [Gu10], Example 7.2. 

\begin{Proposition}
With the same notations as above, we have
$$\supp(\mu_\varphi)=S(E).$$
\end{Proposition}

\textbf{Proof:} See [Gu10], Corollary 7.3. \hfill $\square$

\vspace{0.4cm}

To prove Theorem 1.5 we assume that $\widehat{h}_f$ has no positive lower bound on $K^{nr,v} \setminus \PrePer(f)$. Then there are elements $\{ P_n \}_{n \in \mathbb{N}}$ in $K^{nr,v} \setminus \PrePer(f)$, such that
$$\lim_{n \rightarrow \infty} \widehat{h}_f (P_n) \longrightarrow 0.$$ 

We will show that this contradicts \eqref{support} and Proposition 1.8. In our setting, Yuan's equidistribution Theorem states the following.

\begin{Theorem}[Yuan]
The Galoisorbits of $\{P_n \}$ are equidistributed in $(\mathbb{P}_{\mathbb{K}}^{1})^{an}$. This means:
\begin{eqnarray}
   \mu_f = \lim_{n\rightarrow \infty} \frac{1}{\vert O(P_n) \vert} \sum_{P'\in O(P_n)} \delta_{P'}, \label{Yuan}
  \end{eqnarray}
where $O(P_n )$ denotes the Galoisorbit of $P_n$ and $\delta_{P'}$ the Dirac-measure at $P'$.
\end{Theorem}

\textbf{Proof:} See [Yu], Theorem 3.1. \hfill $\square$

\vspace{0.4cm}

\begin{Korollar}
If there is a sequence $\{ P_n \}_{n\in\mathbb{N}}$ as above, then
$$\supp(\mu_f) \subseteq \val^{-1}(\frac{1}{e_{v\mid p}}\mathbb{Z}\cup \{\pm \infty \}),$$
where $e_{v\mid p}$ is the ramification index of $v$ over $p$. 
\end{Korollar}

\textbf{Proof:} Let $y \in (\mathbb{P}_{\mathbb{K}}^{1})^{an}$ with $\val(y)\notin \frac{1}{e_{v\mid p}}\mathbb{Z}\cup\{\pm \infty\}$. Choose an open neighbourhood $I$ of $\val(y)$, such that $I$ doesn't contain an element of $\frac{1}{e_{v\mid p}}\mathbb{Z}$. The value group of $w$ on $K^{nr,v}$ is $\frac{1}{e_{v\mid p}}\mathbb{Z}$ and $\val$ is continuous, so the open neighbourhood $U_y := \val^{-1}(I)$ of $y$ doesn't contain a rational point of $(\mathbb{P}^{1}_{\mathbb{K}})^{an} (K^{nr,v})$. With \eqref{Yuan} we get
$$\mu_f(U_y)=\lim_{n\rightarrow \infty} \frac{1}{\vert O(P_n) \vert} \sum_{P'\in O(P_n)} \delta_{P'}(U_y)=0.$$
So $y$ is no point of $\supp(\mu_f )$. This proves the Corollary. \hfill $\square$

\vspace{0.4cm}

With Proposition 1.8 and \eqref{support} we conclude
\begin{eqnarray} \pi^{an}(S(E))\subseteq \val^{-1}(\frac{1}{e_{v\mid p}}\mathbb{Z}\cup \{\pm \infty \}) . \label{Widerspruch} \end{eqnarray}

As $\pi^{an}$ and $\val$ are continuous, $\frac{1}{e_{v\mid p}}\mathbb{Z}\cup\{\pm \infty\}$ is discret and $S(E)$ is not, this is very likely to be impossible. But to prove this we need a better understanding of the map $\pi^{an}$.

In rigid geometry, Tate has described the isomorphism between $\nicefrac{\mathbb{K}^*}{q^{\mathbb{Z}}}$ and $E^{an}(\mathbb{K})$. The $x$ and $y$ coordinate in $E^{an}(\mathbb{K})$ of an element $\zeta \in \nicefrac{\mathbb{K}^*}{q^{\mathbb{Z}}}$ are explicitly given by
 $$x(\zeta)=\sum_{n=-\infty}^{\infty} \frac{q^{n} \zeta}{(1-q^n \zeta)²} - 2\sum_{n=1}^{\infty} \frac{nq^n }{(1-q^n )} $$
$$y(\zeta)=\sum_{n=-\infty}^{\infty} \frac{q^{2n} \zeta^2 }{(1-q^n \zeta)^3 } + \sum_{n=1}^{\infty} \frac{nq^n }{(1-q^n )}.$$
For a proof and further information on this isomorphism we refer to [Si99], V.3 and V.4.

Thus $\pi^{an}$ is defined on rational points of $\nicefrac{(\mathbb{G}^{1}_{m})^{an}}{q^\mathbb{Z}}$ by
$$\pi^{an}(\zeta)=\sum_{n=-\infty}^{\infty} \frac{q^{n} \zeta}{(1-q^n \zeta)²} - 2\sum_{n=1}^{\infty} \frac{nq^n }{(1-q^n )}.$$

As a morphism of strict $\mathbb{K}$-affinoid spaces $\pi^{an}$ is induced by a homomorphism $(\pi^{an})^{\sharp} : \mathbb{K}[X] \rightarrow \mathbb{K}\{X,qX^{-1}\}$ of the related $\mathbb{K}$-affinoid algebras (see [Ber], Chapter 2 and [Bo]). With Tate's isomorphism we know
$$(\pi^{an})^{\sharp} (X) = \sum_{n=-\infty}^{\infty} \frac{q^{n} X}{(1-q^n X)²} - 2\sum_{n=1}^{\infty} \frac{nq^n }{(1-q^n )}.$$

Thus for any $f(X) \in \mathbb{K}[X]$ and any $y\in E^{an}$ we have
$$\vert f(X) \vert_{\pi^{an}(y)} = \left\vert f\left(\sum_{n=-\infty}^{\infty} \frac{q^{n} X}{(1-q^n X)²} - 2\sum_{n=1}^{\infty} \frac{nq^n }{(1-q^n )}\right) \right\vert_y .$$

In order to compute $\val(\pi^{an}(0,r))=-\log\vert X \vert_{\pi^{an}(0,r)}$ for an element $(0,r)\in S(E)$ we have to compute $\vert (\pi^{an})^{\sharp}(X)\vert_{(0,r)}$.
It holds $\vert q^n X \vert_{(0,r)} = \vert q^n \vert_w r < 1$ for all $n\geq0$ and hence
$$\left\vert \frac{q^{n} X}{(1-q^n X)²} \right\vert_{(0,r)} = \vert q^{n} X \vert_{(0,r)}$$
for all $n\geq0$. Obviously we have also
$$r = \vert X \vert_{(0,r)} = \vert q^{0} X \vert_{(0,r)} > \vert q^{1} X \vert_{(0,r)} > \cdots$$
leading us to
\begin{eqnarray}
 \left\vert \sum_{n=0}^{\infty} \frac{q^{n} X}{(1-q^n X)²} \right\vert_{(0,r)} = r . \label{Summe2.1}
\end{eqnarray}

For all negative integers $n$, we have $\vert q^{n}X \vert_{(0,r)}>1$, and hence 
$$\left\vert \frac{q^{n} X}{(1-q^n X)²} \right\vert_{(0,r)} = \left\vert \frac{1}{q^n X} \right\vert_{(0,r)} $$
for all $n<0$. With the trivial inequalities
$$ \left\vert \frac{1}{q^{-1} X} \right\vert_{(0,r)} > \left\vert \frac{1}{q^{-2} X} \right\vert_{(0,r)} > \cdots$$
we conclude
\begin{eqnarray}
 \left\vert \sum_{n=1}^{\infty} \frac{q^{-n} X}{(1-q^{-n} X)²} \right\vert_{(0,r)} = \left\vert \frac{1}{q^{-1} X} \right\vert_{(0,r)} = \vert q \vert_w r^{-1}. \label{Summe2.2}
\end{eqnarray}
The equation
\begin{eqnarray}
 \left\vert 2\sum_{n=1}^{\infty} \frac{nq^n }{(1-q^n )} \right\vert_{(0,r)}= \left\vert 2\sum_{n=1}^{\infty} \frac{nq^n }{(1-q^n )} \right\vert_w = \vert 2q \vert_w \label{Summe2.3}
\end{eqnarray}
similarly follows with elementary properties of non-archimedean absolute values. Since $(0,r)$ is an element of the skeleton, we know $\vert q \vert_w < r < 1$. So \eqref{Summe2.1}, \eqref{Summe2.2} and \eqref{Summe2.3} leads us to
\begin{eqnarray}
 \vert X \vert_{\pi^{an}(0,r)} = \left\vert \sum_{n=-\infty}^{\infty} \frac{q^{n} X}{(1-q^n X)²} - 2\sum_{n=1}^{\infty} \frac{nq^n }{(1-q^n )} \right\vert_{(0,r)} \leq \max\{r, \vert q \vert_w r^{-1} \}. \label{Summe2.0}
\end{eqnarray}

If we choose now $(0,r) \in S(E)$ with $1 < r^2 < \vert q \vert_w$ and $\log r \notin \frac{1}{e_{v\mid p}}\mathbb{Z}$, then the value in \eqref{Summe2.0} is equal to $r$. So there is an element with
$$\val(\pi^{an} ((0,r))) = -\log \vert X \vert_{\pi^{an}(0,r)} = - \log r \notin \frac{1}{e_{v\mid p}}\mathbb{Z}.$$
This contradicts \eqref{Widerspruch}. We have shown, that there are no elements $\{P_n \}$ in $K^{nr,v} \setminus \PrePer(f)$ with $\widehat{h}_f (P_n ) \rightarrow 0$. The finiteness of points $a \in K^{nr,v}$ with $\widehat{h}_f (a)=0$ follows with the same proof, when we assume the existence of infinitely many pairwise distinct points $\{P_n \}_{n \in \mathbb{N}}$ with $\widehat{h}_f (P_n )=0$. This is equivalent to the finiteness of $\PrePer(f) \cap K^{nr,v}$ and proves theorem 1.5. \hfill $\square$

\begin{Bemerkung} \rm
The finiteness statement in Theorem 1.5 doesn't hold if we start with an elliptic curve with (potential) good reduction at $v$. In this case the criterion of Néron-Ogg-Shafarevich states, that infinitely many torsion points are unramified over $v$. This means, that there are infinitely many torsion points $P$ in $E(K^{nr,v})$. With \eqref{height} we get for these points $P \in E(K^{nr,v})$ the property $\widehat{h}_f \circ \pi (P) =0$. As $\pi(P)$ is obviously unramified over $v$ and $\deg(\pi)=2$, we have infinitely many elements $a \in K^{nr,v}$ with $\widehat{h}_f (a) =0$.
\end{Bemerkung}

The next proposition shows, that also the Bogomolov Property in Theorem 1.5 in general fails if we start with an elliptic curve with potential good reduction at $v$.

\begin{Proposition}
Let $E$ be an elliptic curve over $K$ with potential good reduction at $v \mid p$. Let $f$ be a Lattès map associated to $E$ and $[m]$. If $p \nmid m$ then there is a finite extension $K'\mid K$ endowed with a non-archimedean place $w \mid v$ and a sequence $\{ q_n \}_{n\in \mathbb{N}_0}$ in $K'^{nr,w}$ with $$\widehat{h}_f (q_n ) \rightarrow 0 .$$
\end{Proposition}

\textbf{Proof:} Choose $K' \mid K$ endowed with a non-archimedean absolute value $w\mid v$ such that $E$ over $K'$ has good reduction at $w$ and there is at least one non-torsion point in $E(K'^{nr,w})$. Let $\mathcal{E}$ be the Néron-model of $E$ over the valuation ring $R_w$ of $K'$. Since $E$ has good reduction at $w$ and $p \nmid m$, we know that the map $[m]:\mathcal{E}\rightarrow\mathcal{E}$ is étale (see [BLR], 7.3 Lemma 2b) and [Si99], Example IV.3.1.4). Especially $[m]$ is unramified.

Now we choose a non-torsion point $Q_0 \in E(K'^{nr,w})$, then $q_0 := \pi (Q_0 )$ is no preperiodic point of $f$. Let $Q_1 \in E(\overline{K'^{nr,w}})$ be a pre-image of $Q_0$ under $[m]$. Using $[m]: \mathcal{E} \rightarrow \mathcal{E}$ is unramified, we get $Q_1 \in E(K'^{nr,w})$ ([BG], Proposition B.3.6).

By successive repetition of this, when we replace $Q_n$ by $Q_{n+1}$ in each step, we get a sequence $\{ Q_n \}_{n\in\mathbb{N}_0}$ in $E(K'^{nr,w})$ with $[m]^n Q_n = Q_0$. We set $q_n := \pi(Q_n )$. By using the Lattès diagram we get $f^n (q_n )=q_0$ for all $n \in \mathbb{N}_0$. As $Q_n$ is in $E(K'^{nr,w})$, $q_n$ is in $K'^{nr,w}$. So we have found a sequence $\{ q_n \}_{n\in\mathbb{N}_0}$ in $K'^{nr,w}$, such that
$$\widehat{h}_f (q_n ) = \left( \frac{1}{m^2} \right)^n \widehat{h}_f (q_0 ) \rightarrow 0$$
(see Definition 1.1). This proves, that the Bogomolov Property cannot hold in this case. \hfill $\square$  


\small Lukas Pottmeyer, Fachbereich Mathematik, Universität Tübingen, D-72076 Tübingen, \texttt{lukas.pottmeyer@uni-tuebingen.de}

\end{document}